\documentclass[11pt]{article}
\usepackage{latexsym, amscd, amsmath,amssymb}

\makeindex
\title{Normal subgroup growth in free class-$2$-nilpotent groups}
\author{Christopher Voll\thanks{Mathematical Institute, Oxford OX1 3LB, England.}}

\newenvironment{acknowledgements}{\bigskip\noindent\textsl{Acknowledgements.}\rm}

\newenvironment{proof}{\bigskip\noindent\textbf{Proof.}\rm}{\hfill$\Box$ \\ \hspace{0.1cm}}

\newtheorem{lemma}{Lemma}
\newtheorem{theorem}{Theorem}

\newtheorem{corollary}{Corollary}

\newtheorem{definition}{Definition}

\long\def\symbolfootnote[#1]#2{\begingroup%
\def\thefootnote{\fnsymbol{footnote}}\footnote[#1]{#2}\endgroup}

\def \Q {\ensuremath{\mathbb{Q}}}

\def \N {\ensuremath{\mathbb{N}}}

\def \Qp {\mathbb{Q}_p}
\def \Z {\mathbb{Z}}
\def \z {\zeta^\triangleleft}
\def \Zp  {\mathbb{Z}_p}

\def \Fp  {\mathbb{F}_p}

\def \LieL {\ensuremath{\mathfrak{L}}}

\def \F23 {\ensuremath{F_{2,3}}}

\def \bfX {{\bf X}}
\def \bfY {{\bf Y}}
\def \bfZ {{\bf Z}}
\def \T {\mathfrak{T}}
\begin{document}
\maketitle
\begin{abstract}
Let $F_{2,d}$ denote the free class-$2$-nilpotent group on~$d$ generators. We
compute the normal zeta functions $\zeta^\triangleleft_{F_{2,d}}(s)$,
prove that they satisfy local functional equations and determine their
abscissae of convergence and pole orders.
\end{abstract}

\section{Introduction}
We fix an integer $d\geq 2$ and write $F_{2,d}$ for the free
nilpotent group of class~$2$ on~$d$ generators. In this paper we compute
the {\sl normal zeta function} of~$F_{2,d}$
$$\zeta^\triangleleft_{F_{2,d}}(s):=\sum_{H\triangleleft_f
  F_{2,d}}|F_{2,d}:H|^{-s}$$
where the sum ranges over all normal subgroups
of~$F_{2,d}$ of finite index and~$s$ is a complex variable. Zeta
functions have been introduced into the study of finitely generated
torsion-free nilpotent ($=:\T$-)groups by Grunewald, Segal and Smith~\cite{GSS/88} as a tool
to study the subgroup growth of these groups. An important invariant
of a $\T$-group~$G$ is the {\sl abscissa of convergence} of its normal
zeta function. Writing
$$s^\triangleleft_{G,n}:=\sum_{i\leq n}a^\triangleleft_{G,i}$$ where
$$a_{G,i}^\triangleleft:=|\{H\triangleleft G|\;|G:H|=i\}|$$
and 
$$\alpha^\triangleleft(G):=\inf\{\alpha\geq0|\;\exists c>0\,\forall n:\;s^\triangleleft_{G,n}<cn^\alpha\},$$
it is well-known that the Dirichlet
series
$$\zeta^\triangleleft_{G}(s)=\sum_{n=1}^\infty a^\triangleleft_{G,n}n^{-s}$$ 
converges on the complex
right half-plane
$$\{s\in\mathbb{C}|\;\mathfrak{Re}(s)>\alpha^\triangleleft(G)\}.$$
As noted in~\cite{GSS/88}, to determine
$\alpha^\triangleleft(G)$ in terms of structural invariants of the
$\T$-group~$F$ is a difficult
problem. In general not much seems to be
known about these numbers. In~\cite{duSG/00} Grunewald and du~Sautoy proved
that they are always rational numbers. They also show that the (normal) zeta function of
a~$\T$-group~$G$ always allows for analytic continuation beyond its abscissa of
convergence---though in general not to the whole complex plane---and that the
continued function has a unique pole of order~$\beta^{\triangleleft}(G)+1$, say, on the line
$\{s\in\mathbb{C}|\;\mathfrak{Re}(s)=\alpha^\triangleleft(G)\}$. Knowledge of
the invariants~$\alpha^\triangleleft(G)$ and $\beta^\triangleleft(G)$ gives us
the following information about the normal subgroup growth of~$G$:

\begin{theorem} [\cite{duSG/00}] \label{duSG theorem}There is a real number~$c$ such that
$$s^\triangleleft_G(n)\sim c\cdot n^{\alpha^\triangleleft(G)}(\log
n)^{\beta^\triangleleft(G)}.$$
\end{theorem}

In Theorem~\ref{thm funeq and abscissa} we will in particular determine these two invariants for
the members of the family~$(F_{2,d})$.

Owing to the nilpotency of the $\T$-group~$G$ the normal zeta
function~$\zeta^\triangleleft_{G}(s)$ decomposes naturally as an Euler product of {\sl local}
normal zeta functions:
$$\zeta^\triangleleft_{G}(s)=\prod_{p\text{
    prime}}\zeta^\triangleleft_{G,p}(s),$$
where $$\zeta^\triangleleft_{G,p}(s):=\sum_{n=0}^\infty
a^\triangleleft_{G,p^n}p^{-ns}.$$
It is a deep result
proved in~\cite{GSS/88} that the local zeta functions of $\T$-groups are rational functions in~$p^{-s}$ . In the case of the free nilpotent groups
$F_{2,d}$ Grunewald, Segal and Smith could prove a stronger result:
\begin{theorem}[\cite{GSS/88}] \label{GSS theorem} The normal zeta
  function $\zeta^\triangleleft_{F_{2,d}}(s)$ is {\rm uniform}, i.e. there exist rational
  functions $W_d^\triangleleft(X,Y)$ of two variables over~$\Q$ such that
  $\zeta^\triangleleft_{F_{2,d},p}(s)=W_d^\triangleleft(p,p^
{-s})$ for all
primes~$p$.
\end{theorem}

In our main result Theorem~\ref{main theorem} we will compute these
rational functions explicitly. As a corollary we get
\begin{theorem} \label{thm funeq and abscissa} We write $F=F_{2,d}$.
\begin{itemize}
\item[a)] The local normal zeta functions satisfy the functional equation
\begin{equation}\label{formula funeq}
\zeta^\triangleleft_{F,p}(s)|_{p\rightarrow p^{-1}}=(-1)^{h(d)}p^{\binom{h(d)}{2}-(h(d)+d)s}\zeta^\triangleleft_{F,p}(s),
\end{equation}
where $h(d)=\binom{d+1}{2}$ is the Hirsch length\symbolfootnote[1]{Recall that the
  Hirsch length of a polycyclic group~$G$ is the number of infinite cyclic
  factors in a decomposition series of~$G$.} of $F_{2,d}$. 
\item[b)]
The
abscissa of convergence of~$\zeta^\triangleleft_{F}(s)$ equals
\begin{equation}
\alpha^\triangleleft(F)=\max\left\{\left.d,\frac{(\binom{d}{2}-j)(d+j)+1}{h(d)-j}\right|\;j\in\{1,\dots,\binom{d}{2}-1\}\right\}.\label{abscissa}
\end{equation}
\item[c)] The zeta function~$\zeta^\triangleleft_{F}(s)$ has a simple pole
  at $s=\alpha^\triangleleft(F)$, i.e.~$\beta^\triangleleft(F)=0$.
\end{itemize}
\end{theorem}

Functional equations of the form~(\ref{formula funeq}) have been found
to hold for all normal zeta
functions of class-$2$-nilpotent groups computed so far. We refer the reader
to~\cite{Voll/03a} for more information on functional equations of zeta
functions which are not finitely uniform. A general group-theoretic explanation of this intriguing symmetry is
still outstanding.

What makes the $\T$-groups $F_{2,d}$ particulary interesting objects
from the point of view of subgroup growth is a connection with
an enumerative problem in finite group theory. In~\cite{duS/02}
M. du~Sautoy showed how the problem of enumerating finite
$p$-groups (on a fixed number~$d$ of genenerators and nilpotency
class~$c$, say) may be reduced to enumerating normal
subgroups in~$F_{c,d}$, the free nilpotent group on~$d$ generators and
class~$c$, {\sl up to the action of the algebraic automorphism group}
of~$F_{c,d}$. The computation of the normal zeta
functions~$\zeta^\triangleleft_{F_{2,d},p}(s)$ may thus be considered as
a step towards the determination of the Dirichlet series
enumerating class-$2$-nilpotent $d$-generator $p$-groups. We refer the
reader to~\cite{duS/02} for details.

\bigskip
The proof of Theorem~\ref{GSS theorem} given in~\cite{GSS/88} reduces to showing that the
coefficients of the local zeta functions~$\zeta^\triangleleft_{F_{2,d},p}(s)$ are polynomials in the
prime~$p$. It might raise the impression that determining the local
zeta functions explicitly would be at least as hard as writing down
formulae for the {\sl Hall polynomials}~$g_{\mu,\nu}^\lambda$. Given
partitions $\lambda, \mu,\nu$, these polynomials in~$p$ count
the number of subgroups of an
abelian $p$-group of type~$\lambda$ which have type~$\mu$ and
  cotype~$\nu$. No simple formula seems to be known for these
polynomials. The computations leading to our main Theorem~\ref{main
  theorem} are based on the observation that it
suffices to know explicitly the polynomial
$$g_{\mu}^\lambda:=\sum_{\nu\leq\lambda}g_{\mu,\nu}^\lambda,$$
for which a very simple (if apparently little-known) formula has been known for at
least a century (see Lemma~\ref{lemma subgroups of abelian groups} in Section~\ref{section proof of
  main theorem}). For our computations we use the Cartan decomposition of certain
lattices in the Lie ring associated to~$F_{2,d}$, rather than the
machinery of $p$-adic integration.

Our main Theorem~\ref{main theorem} presents the local zeta
functions~$\z_{F_{2,d},p}(s)$ as a finite sum of rational functions in~$p$
and~$p^{-s}$. We will state it in Section~\ref{section notation and
  main result}
and prove it in Section~\ref{section proof of main theorem}. The proof of Theorem~\ref{thm funeq and
  abscissa}, part $a)$,  reduces to showing that each summand satisfies a
functional equation like~(\ref{formula funeq}). The proof of part~$b)$ depends on
identifying a `dominating' summand that determines the
abscissa of convergence of the global zeta function~$\z_{F_{2,d}}(s)$. 
Part~$c)$ is a consequence of the fact that for any given~$d$ the maximum
defining~$\alpha^\triangleleft(F_{2,d})$ in~(\ref{abscissa}) is
obtained at a {\sl unique} value of~$j$. This in turn follows from the
non-existence of certain integral points on a particular elliptic curve.

The zeta functions $\z_{F_{2,d}}(s)$ for $d\in\{2,3\}$ were computed
in~\cite{GSS/88}. The zeta function $\z_{F_{2,4}}(s)$ was first
computed by Paajanen~\cite{Paajanen/04}, extending methods introduced
in~\cite{Voll/03a}. She also observed that, for any integer~$d$, the rational number
specified in~(\ref{abscissa}) is a good
lower bound for~$\alpha^\triangleleft(F_{2,d})$ and, {\sl mutatis mutandis},
for the abscissa
of convergence of the normal zeta function of any class-$2$-nilpotent group. Note that $F_{2,5}$ is the first member of the
family~$(F_{2,d})$ for which the abscissa of convergence $\alpha^\triangleleft$ is
not an integer. The best estimates so far for the
numbers~$\alpha^\triangleleft(F_{2,d})$ were derived in~\cite{GSS/88}
where the inequalities
$$(d^3-d^2+2)/4d\leq
\alpha^\triangleleft(F_{2,d})\leq\max\{d,(d-1)(d+1)/2\}$$
were proved.

\section{Notation and statement of the main result}\label{section notation and
  main result}

Throughout this paper we will use the following notation:
\bigskip

\begin{tabular}{ll}
$\N$&$\{1,2,\dots\}$\\
$n'$&$\binom{n}{2}=n(n-1)/2,\; n\in \N$\\
$[n]$&$\{1,\dots,n\},\; n\in \N$\\
$I_0$&$I\cup\{0\},\;I\subseteq\N$\\
$I=\{i_1,\dots,i_m\}_{<}$&$I\subseteq [n],\;1\leq i_1<\dots<i_m\leq
n$\\
$\binom{a}{b}_p$&$\prod_{i=0}^{b-1}(p^a-p^i)/(p^b-p^i),\;a,b\in\N,
  a\geq b, p \text{ prime,}$\\
&the number of subspaces of dimension $b$ in $\Fp^a$\\
$\binom{a}{b}_{p^{-1}}$&$\left.\binom{a}{b}_p\right|_{p\rightarrow p^{-1}}$\\
$b_{n,I}(p)$&$\binom{n}{i_m}_p\binom{i_m}{i_{m-1}}_p...\binom{i_2}{i_1}p,$\\
& the number of flags of dimension type~$I$ in~$\Fp^n$\\
$\Zp$&  the ring of $p$-adic integers\\
$\Qp$&  the field of $p$-adic numbers\\
$\zeta(s)$ & $\sum_{n=1}^\infty n^{-s}$, the Riemann zeta function \\
$\zeta_p(s)$ & $\sum_{n=0}^\infty p^{-ns}$, the local Riemann zeta function at
the prime~$p$.
\end{tabular}
\bigskip

The following definitions will allow us to state our main
Theorem~\ref{main theorem}. Motivations for
these definitions will be provided in Section~\ref{section proof of
  main theorem} where this result will be proved.

Given  $n\in\N$, a prime~$p$ and a vector of independent variables
$\bfZ=(Z_1,\dots,Z_{n-1})$ we set
$$ F_n(p,\bfZ):=\sum_{I\subseteq[n-1]}b_{n,I}(p^{-1})\prod_{i\in
  I}\frac{Z_{i}}{1-Z_{i}}.$$
Given integers $0\leq a \leq b\leq n$ we will write
$F_{b-a}(p,\bfZ_{a}^{b})$ for $F_{b-a}(p,(Z_{{a}+1},\dots,Z_{b-1}))$
and $F_{b-a}(p,\bf{Z^{-1}}_{a}^{b})$ for $F_{b-a}(p,(Z^{-1}_{{a}+1},\dots,Z^{-1}_{b-1}))$.

We set
\begin{eqnarray}
\phi:[d-1]_0&\rightarrow&[d']\label{definition phi}\\
      i&\mapsto &id-(i+1)'. \nonumber
\end{eqnarray}
Given subsets $I=\{i_1,\dots,i_h\}_<\subseteq[d-1]_0$ and
    $J=\{j_1,\dots,j_h\}_<\subseteq[d']$ of equal cardinality we will write $\phi(I)\leq
    J$ if $\forall r\in [h]:\;\phi(i_r)\leq j_r$.

A choice of subsets $I=\{i_1,\dots,i_h\}_<\subseteq A\subseteq\N$ and
$J=\{t_1,\dots,t_h\}_<\subseteq B\subseteq\N$ allows us to define a total order $\prec(=\prec_{I,J})$ on the
disjoint union $A\uplus B$ as follows. For $x,y\in
A\uplus B, a\not=b$ we set

\begin{eqnarray*}
 x\prec y &\text{ if}&\left\{\begin{array}{lll}x<y&\text{ if }&(x,y \in A)\vee (x,y \in
  B),\\\exists r\in[h]:\;x\leq i_r \wedge y\geq j_r&\text{ if }&(x\in
  A)\wedge(y\in  B),\\\forall r\in[h]:\;j_r\leq x\Rightarrow i_r< y&\text{ if }& (y\in A)\wedge(x\in  B),\end{array}\right.
\end{eqnarray*}
where~$<$ refers to the natural order on~$\N$.
If, for example, $A=B=\N$ the list of elements of $(\N\uplus
\N,\prec_{I,J})$ in ascending order starts like this:
$$1_B\prec\dots\prec (j_1-1)_B\prec 1_A\prec\dots\prec (i_1-1)_A\prec
(j_1)_B\prec\dots\prec (j_2-1)B\prec (i_1)_A\prec\dots$$

We can now state our main theorem. Let $F=F_{2,d}$.
\begin{theorem}\label{main theorem} For all primes $p$
$$\zeta^\triangleleft_{F,p}(s)=\zeta_{\Zp^d}(s)\zeta_p((d'+d)s-dd')A(p,p^{-s}),$$ where
\begin{equation}
A(p,p^{-s}):=\sum_{\substack{I\subseteq[d-2], J\subseteq [d'-1]\\ |I|=|J|,
    \phi(I)\leq J}}A_{I,J}(p,p^{-s})\label{A}
\end{equation}
and, if $I=\{i_1,\dots,i_h\}_<$, $J=\{j_1,\dots,j_h\}_<$, 
\begin{eqnarray}
\lefteqn{A_{I,J}(p,p^{-s}):=F_{d'-{j_h}}(p,\bfX_{j_h}^{d'})\zeta_{\Zp^{i_h}}^{-1}(s)}\cdot\label{A_IJ}\\
&&\prod_{r=1}^{h}
\binom{j_{r+1}-\phi(i_r)}{j_r-\phi(i_r)}_{p^{-1}}\binom{d-i_{r-1}}{d-i_{r}}_{p^{-1}} \frac{Y_{i_r}}{(1-Y_{i_r})(1-Y'_{i_r})}F_{j_r-j_{r-1}}(p,\bfX_{j_{r-1}}^{j_r})F_{i_r-i_{r-1}}(p,\bfY_{i_{r-1}}^{i_r})\nonumber
\end{eqnarray}
with numerical data $(X_{j})$, $(Y_{i})$, $(Y'_{i_r})$ defined by
\begin{eqnarray}
X_{j}&:=&p^{i(j)(d-i(j))+(d'-j)(d+j-\phi(i(j)))-(d-i(j)+d'-j)s},\quad\;j\in [d'-1] \label{num
  data}\\
Y_{i}&:=&p^{i(d-i)+(d'-j(i))(d+j(i)-\phi(i))-(d-i+d'-j(i))s},\quad\;i\in [d-2]\nonumber\\
Y'_{i_r}&:=&p^{i_{r-1}(d-i_{r-1})+(d'-j(i_{r}))(d+j(i_{r})-\phi(i_{r-1}))-(d-i_{r-1}+d'-j(i_{r}))s},\quad\;r\in[h]\nonumber
\end{eqnarray}
where 
 \begin{eqnarray*}
j(i)&:=&\min\{j\in J\cup\{d'\}|\;\phi(i)\prec_{\phi(I),J} j\},\\
 i(j)&:=&\max\{i\in
I\cup\{0\}|\;\phi(i)\prec_{\phi(I),J} j\}
\end{eqnarray*}
and $j_0:=i_0:=0\prec_{\phi(I),J} \phi([d-2])\uplus [d'-1]\prec_{\phi(I),J} j_{h+1}:=d'$.
\end{theorem}
Note that $i(j_r)=i_r$ and $j(i_r)=j_r$ and thus $X_{j_r}=Y_{i_r}$ and that,
owing to our assumption on~$I$ and~$J$, we have~$\forall
 j\in[d'-1]\;\phi(i(j))\leq j$. Recall that
$$\zeta_{\Zp^d}(s)=\prod_{i=0}^{d-1}\zeta_p(s-i).$$

\section{Proof of Theorem~\ref{main theorem}}\label{section proof of
  main theorem}
Let $n\in\N$. Two lattices $\Lambda,\Lambda'\subseteq\Zp^n$ are said to be
{\sl homothetic} if there is a non-zero constant~$c\in\Qp$ such that
$c\Lambda= \Lambda'$. A lattice~$\Lambda\subseteq\Zp^n$ is called {\sl maximal} (in
its homothety class~$[\Lambda]$) if~$p^{-1}\Lambda\not\subseteq\Zp^n$.
It is said to be of {\sl
   type}~$\nu(\Lambda)=(I,{\bf r}_{I_0})$ if 
\begin{equation}
I=\{i_1,\dots,i_{l}\}_<\subseteq[{n}],\;{\bf r}_{I_0}=(r_0,r_{i_1},\dots,r_{i_{l}}
),\label{I}
\end{equation}
and $r_{0}\in\N_0$, $r_{i_j}\in\mathbb{N}\mbox{ for }i_j\in I$ and $\Lambda$ has elementary divisors
\begin{equation}
\left(\underbrace{p^{\sum_{i\in I_0}r_{i}},\dots,{p^{\sum_{i\in I_0}r_{i}}}}_{n-i_{l}},\dots,\underbrace{p^{r_{0}+r_{i_1}},\dots,p^{{r_{0}+r_{i_1}}}}_{i_2-i_1},\underbrace{p^{r_{0}},\dots,p^{{r_{0}}}}_{i_1}\right)=:(p^{\nu}).\label{eldiv type}
\end{equation}

By slight abuse of notation we may  say that a  lattice is of
type~$I$ if it is of type~$(I,{\bf r}_{I_0})$ for some
vector~${\bf r}_{I_0}$ and that the {\sl homothety class}~$[\Lambda]$ has
type~$I$ if some representative has type~$I$, in which case we
write~$\nu([\Lambda])=I$.

Now let $F=F_{2,d}$ and, for a prime~$p$, define the $\Zp$-Lie algebra
(with Lie-bracket induced from taking commutators in $F$) $$\LieL_p:=(F/Z(F)\oplus
Z(F))\otimes_{\Z}\Zp,$$
where $Z(F)$ is the centre of the group $F$. Denote by $\LieL_p'=Z(\LieL_p)$
the centre of~$\LieL_p$. As in the proof of Theorem~\ref{GSS theorem} in~\cite{GSS/88} we note that
if~$\Lambda\subseteq \LieL_p$ and
$\overline{\Lambda}:=(\Lambda+\LieL'_p)/{\LieL'_p}\subseteq \LieL_p/{\LieL'_p}\cong\Zp^d$ is of
type $$(S,{\bf q}_{S_0}),\; S=\{s_1,\dots,s_m\}_<,\;{\bf
  q}_{S_0}=(q_0,q_{s_1},\dots,q_{s_{m}})$$ then
$\phi(\Lambda):=[\Lambda,\LieL_p]\subseteq \LieL'_p$ is of
type\symbolfootnote[1]{Note that for this statement to make sense we had to
  allow $n\in I$ in~(\ref{I}), as $\phi(d-1)=d'$.}
$$(\phi(S),{\bf q}_{\phi(S)_0})$$
with~$\phi$ as defined
by~(\ref{definition phi}) in Section~\ref{section notation and
  main result}. The proof concluded by
showing that the number of (WLOG maximal) lattices $M\subseteq \LieL'_p$
containing $\phi(\Lambda)$ and of type~$(T,{\bf r}_T)$, $T\subseteq[d'-1]$,
is a polynomial in~$p$. Giving an expression for this polynomial in
terms of Hall polynomials rather obscures its simple shape, given in
the following
Lemma~\ref{lemma subgroups of abelian groups}. Given two
partitions $\mu\subseteq\lambda$ we denote by~$\alpha_{\lambda}(\mu;p)$
the number of subgroups of type~$\mu$ in a finite abelian $p$-group of
type~$\lambda$ (where an abelian~$p$-group~$G$ is said to be of
type~$\lambda=(\lambda_1,\dots,\lambda_n)$,
$\lambda_1<\dots<\lambda_n$ if $G\cong\Z/(p^{\lambda_1})\times\dots\times\Z/(p^{\lambda_n})$).

\begin{lemma} [1.4.1 in~\cite{Butler/94}] \label{lemma subgroups of abelian groups}
  For any two partitions $\mu\subseteq\lambda$,
$$\alpha_{\lambda}(\mu;p)=\prod_{j\geq1}p^{\mu'_{j}(\lambda'_j-\mu'_j)}\binom{\lambda'_j-\mu'_{j+1}}{\lambda'_j-\mu'_{j}}_{p^{-1}}$$
where $\lambda',\mu'$ are the conjugate partitions of $\lambda,\mu$, respectively.
\end{lemma}

For our computations it will be useful to associate to a
pair~$(\Lambda,M)$ of lattices with~$\phi(\Lambda)\subseteq M$ an invariant measuring the `overlap' of their types.

\begin{definition} \label{definition I,J}
Given two lattices $\Lambda \subseteq \LieL_p$ and $M\subseteq
\LieL'_p$
such that~$\phi(\Lambda)\subseteq M$, $M$ maximal and
\begin{alignat*}{2}
\overline{\Lambda}&\text{ is of type }(S,{\bf
  q}_{S_0}),\quad&S&=\{s_1,\dots,s_m\}\subseteq[d-1],\\
M&\text{ is of type }(T,{\bf r}_T),&T&=\{t_1,\dots,t_n\}\subseteq[d'-1]
\end{alignat*}
 and with
$$Q_{k}:=q_0+\sum_{k\geq j\in[m]}q_{s_j},\; k\in[m]_0, \quad\text{
  and }\quad R_{l}:=\sum_{l\geq i\in [n]}r_{t_i},\; l\in[n]_0,$$ we set
\begin{eqnarray*}
I=I(\Lambda,M)&:=&\{s_k\in S|\;\exists l\in[n]:\;R_{{l-1}}\leq
Q_{{k}}<R_{l}\}\\
J=J(\Lambda,M)&:=&\{t_l\in T|\;\exists k\in[m]:\;R_{{l-1}}\leq
Q_{{k}}<R_{l}\}.
\end{eqnarray*}
\end{definition}
Note that always $|I|=|J|$, $I\subseteq[d-2]$, $J\subseteq[d'-1]$ and $\phi(I)\leq J$. Conversely, every such pair~$(I,J)$ occurs
in this way. We may thus write
\begin{eqnarray}
\zeta^\triangleleft_{F,p}(s)=\zeta_p((d'+d)s-dd')B(p,p^{-s}),\label{zeta
  B}
\end{eqnarray}
 where
\begin{equation*}
B(p,p^{-s}):=\sum_{\substack{I\subseteq[d-2], J\subseteq[d'-1]\\ |I|=|J|,
    \phi(I)\leq J}}B_{I,J}(p,p^{-s})\label{B}
\end{equation*}
and
\begin{equation}
B_{I,J}(p,p^{-s}):=\sum_{\substack{\Lambda\subseteq \LieL_p,
    \phi(\Lambda)\subseteq M\subseteq \LieL'_p\\ M\text{
    maximal}\\I(\Lambda,M)=I,
    J(\Lambda,M)=J}}|\LieL_p/\LieL'_p:\overline{\Lambda}|^{-s}|\LieL'_p:M|^{d-s}.
    \label{zeta B def}
\end{equation}
Note that the factor $\zeta_p((d'+d)s-dd')$ in~(\ref{zeta B}) allowed us
to sum just over {\sl maximal} lattices in~(\ref{zeta B def}).
To apply Lemma~\ref{lemma subgroups of abelian groups} we note that
the partition~$\mu'$ of~$\log_p|\LieL_p':M|$ associated to the maximal
lattice~$M$ equals
$$\left((d'-t_1)^{r_{t_1}},\dots,(d'-t_n)^{r_{t_n}}\right),$$
while the partition~$\lambda'$ of~$\log_p|\LieL_p':\phi(\Lambda)|$
associated to~$\phi(\Lambda)$ equals
$$\left((d')^{q_0},(d'-\phi(s_1))^{q_{s_1}},\dots,(d'-\phi(s_m))^{q_{s_m}}\right).$$
Thus, if $I=\{i_1,\dots,i_h\}_<\subseteq[d-2]$ and $J=\{j_1,\dots,j_h\}_<\subseteq[d'-1]$, 
\begin{eqnarray*}\lefteqn{B_{I,J}(p,p^{-s})=}\\
&&\sum_{B\subset[j_1]}b_{j_1,B}(p^{-1})\sum_{{\bf
    r}_B>0}p^{\sum_{b\in
    B}r_b\left((d+b)(d'-b)-s(d+d'-b)\right)}\sum_{e=0}^\infty
    p^{e\left((d'-j_1)(j_1+d)-s(d_d'-j_1)\right)}\\
&&\sum_{A\subset[i_1]}b_{i_1,A}(p^{-1})\sum_{{\bf q}_A>0}p^{\sum_{a\in
    A}q_a\left(a(d-a)+(d'-j_1)(d'+j_1-\phi(a)) -s(d+d'-a-j_1)\right)}\\
&&\binom{d}{d-i_1}_{p^{-1}}\binom{j_2-\phi(i_1)}{j_1-\phi(i_1)}_{p^{-1}}\sum_{f=1}^\infty
    p^{f\left(i_1(d-i_1)+(d'-j_1)(d+j_1-\phi(i_1))-s(d+d'-j_1)\right)}\\
&&\left.\dots.\right.\\
&&\sum_{B\subset[j_h]\setminus[j_{h-1}]}b_{j_h-j_{h-1},B-j_{h-1}}(p^{-1})\sum_{{\bf
    r}_B>0}p^{\sum_{b\in
    B}r_b\left(i_{h-1}(d-i_{h-1})+(d+b-\phi(i_{h-1}))(d'-b)-s(d+d'-b-i_{h-1})\right)}\\
&&\sum_{e=0}^\infty
    p^{e\left(i_{h-1}(d-i_{h-1})+(d'-j_h)(j_h-j_{h-1}+d)-s(d_d'-\phi(i_{h-1})-j_1)\right)}\\
&&\sum_{A\subset[i_h]\setminus[i_{h-1}]}b_{i_h-i_{h-1},A-i_{h-1}}(p^{-1})\sum_{{\bf q}_A>0}p^{\sum_{a\in
    A}q_a\left(a(d-a)+(d'-j_h)(d'+j_h-\phi(a)) -s(d+d'-a-j_h)\right)}\\
&&\binom{d-i_{h-1}}{d-i_h}_{p^{-1}}\binom{d-\phi(i_h)}{j_h-\phi(i_h)}_{p^{-1}}\sum_{f=1}^\infty
    p^{f\left(i_h(d-i_h)+(d'-j_h)(d+j_h-\phi(i_h))-s(d+d'-j_h-i_{h-1})\right)}\\
&&\sum_{B\subset[d']\setminus[j_{h}]}b_{d'-j_{h},B-j_h}(p^{-1})\sum_{{\bf
    r}_B>0}p^{\sum_{b\in
    B}r_b\left(i_{h}(d-i_{h})+(d+b-\phi(i_{h}))(d'-b)-s(d+d'-b-i_{h})\right)}\\
&&\underbrace{\sum_{e=0}^\infty
    p^{e(i_h(d-i_h)-s(d-i_h))}\sum_{A\subset[d]\setminus[i_h]}b_{d-i_h,A-i_h}(p^{-s})\sum_{{\bf
    q}_A>0}p^{\sum_{a\in A}q_a\left(a(d-a)-s(d-a)\right)}}_{=\zeta_p(s-i_h)\dots\zeta_p(s-d+1)=\zeta_{\Zp^d}(s)/\zeta_{\Zp^{i_h}}(s)}.
\end{eqnarray*}
(Here ${\bf q}_A>0$ stands for $(q_a)_{a\in A}\in\N^{|A|}$. Similarly
for ${\bf r}_B>0$.)
One checks easily that the theorem follows with 
$$A_{I,J}(p,p^{-s}):=B_{I,J}(p,p^{-s})/\zeta_{\Zp^d}(s).$$

\section{Proof of Theorem~\ref{thm funeq and abscissa}}\label{section proof of thm funeq and abscissa}
\subsection{Part $a)$} \label{section proof of funeq}
To prove the functional equation~(\ref{formula funeq}) it clearly
suffices to show that for every summand in~(\ref{A})
\begin{equation} \label{funeq A_IJ}
\forall\;
I\subseteq[d-2],J\subseteq[d'-1]:\;\left.A_{I,J}(p,p^{-s})\right|_{p\rightarrow
  p^{-1}}=(-1)^{d'-1}p^{\binom{d'}{2}}A_{I,J}(p,p^{-s}).
\end{equation}
Using the observation that $\binom{a}{b}_{p^{-1}}=p^{b(a-b)}\binom{a}{b}_{p}$ it is
not hard to verify that
\begin{eqnarray} \label{G_r plus binos}
\lefteqn{\log_p\left(\prod_{r=1}^h\frac{\binom{j_{r+1}-\phi(i_r)}{j_r-\phi(i_r)}_{p}\binom{d-i_{r-1}}{d-i_{r}}_{p}
  \frac{Y_{i_r}^{-1}}{(1-Y_{i_r}^{-1})(1-{Y'}_{i_r}^{-1})}}{\binom{j_{r+1}-\phi(i_r)}{j_r-\phi(i_r)}_{p^{-1}}\binom{d-i_{r-1}}{d-i_{r}}_{p^{-1}}
    \frac{Y_{i_r}}{(1-Y_{i_r})(1-Y'_{i_r})}}\right)=}\\
&&
    \left(\sum_{r=1}^hi_{r-1}(i_r-i_{r-1})+j_r(j_{r+1}-j_r)\right)-i_hs. \nonumber
\end{eqnarray}
(Recall that we had set $j_0=i_0=0$ and $j_{h+1}=d'$.) As
$$\zeta_{\Zp^n}(s)=\prod_{i=0}^{n-1}\zeta(s-i)$$ for $n\in\N$ we
have of course
$$\log_p\left((-1)^{i_h}\frac{\zeta_{\Zp^{i_h}(s)}}{\zeta_{\Zp^{i_h}(s)}|_{p\rightarrow p^{-1}}}\right)=-\binom{i_h}{2}+i_hs.$$
Repeatedly applying the equation (proved, for example, in~\cite{Voll/03a})
$$ F_n(p^{-1},\bfZ^{-1})=(-1)^{n-1}p^{\binom{n}{2}} F_n(p,\bfZ)$$
we get
\begin{eqnarray*}
\log_p\left((-1)^{i_h-h}\prod_{r=1}^h
  \frac{F_{i_r-i_{r-1}}(p^{-1},{{\bfY}^{-1}}_{i_{r-1}}^{i_r})}{F_{i_r-i_{r-1}}(p,\bfY_{i_{r-1}}^{i_r})}\right)&=&\sum_{r=1}^h\binom{i_r-i_{r-1}}{2},\\
\log_p\left((-1)^{d'-1-h}\prod_{r=0}^h
  \frac{F_{j_r-j_{r-1}}(p^{-1},{{\bfX}^{-1}}_{j_{r-1}}^{j_r})}{F_{j_r-j_{r-1}}(p,\bfX_{j_{r-1}}^{j_r})}\right)&=&\sum_{r=0}^h\binom{j_{r+1}-j_{r}}{2}.
\end{eqnarray*}
But
\begin{eqnarray*}
\sum_{r=1}^h\binom{i_r-i_{r-1}}{2}+i_{r-1}(i_r-i_{r-1})&=&\sum_{r=1}^h\binom{i_r}{2}-\binom{i_{r-1}}{2}=\binom{i_h}{2},\\
\sum_{r=0}^h\binom{j_{r+1}-j_{r}}{2}+j_{r}(j_{r+1}-j_{r})&=&\binom{j_{h+1}}{2}=\binom{d'}{2}.
\end{eqnarray*}
Equation~(\ref{funeq A_IJ}) follows. 

\subsection{Part $b)$}  \label{section proof of alpha}
The simplest of the rational functions $A_{I,J}(p,p^{-s})$ defined
by~(\ref{A_IJ}) is 
$$A_0(p,p^{-s}):=A_{\emptyset,\emptyset}(p,p^{-s})=F_{d'}(p,\bfX),$$
with $\bfX=(X_1,\dots,X_{d'-1})$ and
$$X_j=p^{(d'-j)(d+j)-(d+d'-j)s}.$$
Part $b)$ of Theorem~\ref{thm funeq and abscissa} will be proved if we
can show that this summand determines the abscissa of convergence
of~$\z_{F}(s)$. This will follow once we show that the rational function~$A$
in Theorem~\ref{main theorem} satisfies the conditions of the following lemma.

\begin{lemma}
Let $A\in\Q(X,Y)$ be a rational function in two variables over~$\Q$ with
the property that
\begin{enumerate}
\item[a)] $A=\sum_{i=0}^n A_i,\;A_i=\frac{P_i}{Q_i},\; P_i,Q_i\in\Q[X,Y]$.
\item[b)] $\forall
  i\in[n]_0:\;Q_i(X,Y)=\prod_{j=1}^{r_i}(1-X^{a_{ij}}Y^{b_{ij}})$ for
 integers $a_{ij}\in\N_0,b_{ij},r_i\in\N$.
\item[c)] $P_0(X,Y)=1+\sum_{j=1}^{s_0}X^{c_{0j}}Y^{d_{0j}}$ for
   integers $c_{0j}\in\N_0,d_{0j},s_0\in\N$.
\item[d)] $\forall
  i\in[n]:\;P_i(X,Y)=\sum_{j=1}^{s_j}X^{c_{ij}}Y^{d_{ij}} \;$ for integers $c_{ij}\in\N_0,d_{ij},s_i\in\N$.
\item[e)] $\max_{\substack{i\in[n]_0\\j\in[s_i]}}\left\{\frac{c_{ij}+1}{d_{ij}}\right\}<\max_{\substack{i\in[n]_0\\j\in [r_i]}}\left\{\frac{a_{ij}+1}{b_{ij}}\right\}=\max_{j\in [r_0]}\left\{\frac{a_{ij}+1}{b_{ij}}\right\}$.
\end{enumerate}
Then $\max_{j\in [r_0]}\left\{\frac{a_{ij}+1}{b_{ij}}\right\}$ is the abscissa of convergence of the Euler product $\prod_{p
  \text{ prime}}A(p,p^{-s})$.
\end{lemma}
\begin{proof}
Writing
$$\prod_{p
  \text{ prime}}A(p,p^{-s})=\prod_{\substack{i\in[n]_0\\j\in
    [r_i]}}\zeta(b_{ij}s-a_{ij})\prod_{p \text{ prime}}N(p,p^{-s}),$$
say, it is enough to prove that $\prod_{p \text{ prime}}N(p,p^{-s})$
converges at $s=\max_{j\in
  [r_0]}\left\{\frac{a_{ij}+1}{b_{ij}}\right\}$, the right-most pole of the
product of Riemann zeta functions on the right hand side. But this
follows from $e)$ and repeated applications of $\frac{a_1+a_2}{{b_1+b_2}}\leq\max\left\{\frac{a_1}{b_1},\frac{a_2}{b_2}\right\}$.
\end{proof}

One checks immediately that the rational function $A$ in
Theorem~\ref{main theorem} satisfies conditions $a)$ to $d)$
with any labelling of the finitely many $A_{I,J}$, $(I,J)\not=(\emptyset,\emptyset)$,
as $A_1,\dots,A_n$. To check condition~$e)$ we need to investigate the
numerical data~(\ref{num data}).  The key observation is that the
complicated dependence of~$X_j$, $Y_i$ and~$Y'_{i_r}$ on the total
order~$\prec_{\phi(I),J}$ need not concern us here: All we need is the
following lemma. The sketch of its trivial proof is only included for
completeness (and certainly not more enlightening than a Computer Algebra plot).

\begin{lemma}
The restriction of the rational function
\begin{eqnarray*}
f_d:[0,d-2]\times[1,d'-1]&\rightarrow&\mathbb{R}\\
(i,j)&\mapsto&\frac{i(d-i)+(d'-j)(d+j-\phi(i))+1}{d+d'-j-i}
\end{eqnarray*}
to $\N_0^2\cap[0,d-2]\times[1,d'-1]$ attains its maximum on the line
$i=0$. 
\end{lemma}
\begin{proof} We will first show that the function $f_d$ attains its
  maximum on the line~$i=0$. Then it will suffice to show that
  $\max\left\{f_d(i,j)|i=0,j\in[d'-1]\right\}>\left.f_d\right|_{[1,d-2]\times[1,d'-1]}$.

First note that, for all~$i\in[0,d-2]$, the restriction
$f_d|_{\left\{i\right\}\times[1,d'-1]}$ attains its unique maximum
  $$f_d(i,j^{(i)})=((i-d)^2-3(i-d))/2+d-\sqrt{-2(i-d)^3+6(i-d)^2-4}$$
  at $$j=j^{(i)}=(d^2+d)/2-i-\sqrt{-2(i-d)^3+6(i-d)^2-4}/2.$$
We want to show that
$$f_d(0,j^{(0)})>f_d(i,j^{(i)})\;\forall i\in]0,d-2].$$
Write $$f_d(i,j^{(i)})=g_d(X),\; X:=i-d,$$
and observe that that $g_d(X)$ has a unique minimum in the
interval~$[-d,-2]$. In fact $g_d(X)'=0$ has a unique solution $X_0$
in~$[-d,-2]$ and one checks easily that
$g_d(X_0)<g_d(-2)<g_d(-d)=f_d(0,j^{(0)})$. Thus $f_d$ attains its
  maximum on the line~$i=0$.                                      

To finish the proof of the lemma it suffices to show that the values at
the integers close to~$j^{(0)}$ are still larger than
$f_d(1,j^{(1)})$. This follows from the equalities
$$f_d(0,j^{(0)})-f_d(1,j^{(1)})=d+1-\sqrt{-4+6d^2+2d^3}+\sqrt{2d^3-6d}$$
and
$$f_d(0,j^{(0)})-f_d(0,j^{(0)}+1)=\frac{2}{\sqrt{-4+6d+2d^3}-2}.$$
\end{proof}

\subsection{Part $c)$}\label{section proof of beta=0}

It clearly suffices to show that there is no integer $j\in\{1,\dots,d'-2\}$ such
that 
\begin{equation} \label{equality}
f_d(0,j)=f_d(0,j+1).
\end{equation}
Indeed, equality~(\ref{equality}) implies 
$$j=\frac{d^2+d-1-\sqrt{2d^3+6d^2-3}}{2}.$$
This fraction is integral only if $2d^3+6d^2-3$ is the square of an odd
integer, i.e. if
\begin{equation}
d^3+3d^2=2(a^2+a+1)\label{elliptic curve}
\end{equation}
has a solution in positive integers. By a Theorem of Siegel there are at most
finitely many such solutions. But, in fact, by coordinate transformations 
$a=-1/2+y/4$, $
d=-1-x/2$
we may write (\ref{elliptic curve}) in Weierstrass form
$$y^2=x^3-12x+4.$$
Neither of the two integral solutions $(x,y)=(0,2)$ and
$(x,y)=(9,25)$ to this equation give us a positive 
integral solution to~(\ref{elliptic curve}).

\begin{acknowledgements} We should like to thank the
  UK's Engineering and Physical Sciences Research Council~(EPSRC) for
  their support in the form of a Postdoctoral Fellowship, and the
  Heinrich-Heine-Universit\"{a}t D\"{u}sseldorf for their hospitality during
  the writing of this paper. We thank
  Fritz Grunewald, Benjamin Klopsch and Odile Sauzet for helpful and inspiring
  conversations. The computer algebra packets Maple and Magma were used in
  Sections~\ref{section proof of alpha} and~\ref{section proof of beta=0}.
\end{acknowledgements}

\bibliographystyle{amsplain}
\bibliography{thebibliography}

\providecommand{\bysame}{\leavevmode\hbox to3em{\hrulefill}\thinspace}
\providecommand{\MR}{\relax\ifhmode\unskip\space\fi MR }
\providecommand{\MRhref}[2]{%
  \href{http://www.ams.org/mathscinet-getitem?mr=#1}{#2}
}
\providecommand{\href}[2]{#2}
\begin{thebibliography}{1}

\bibitem{Butler/94}
L.M. Butler, \emph{Subgroup lattices and symmetric functions}, Mem. Amer. Math.
  Soc. \textbf{112} (1994).

\bibitem{duS/02}
M.P.F. du~Sautoy, \emph{Counting $p$-groups and nilpotent groups}, Publ. Math.
  I.H.E.S. \textbf{92} (2000).

\bibitem{duSG/00}
F.J. Grunewald and M.P.F. du~Sautoy, \emph{Analytic properties of zeta
  functions and subgroup growth}, Ann. Math. \textbf{152} (2000), 793--833.

\bibitem{GSS/88}
F.J. Grunewald, D.~Segal, and G.C. Smith, \emph{Subgroups of finite index in
  nilpotent groups}, Invent. Math. \textbf{93} (1988), 185--223.

\bibitem{Paajanen/04}
P.M. Paajanen, \emph{{The normal zeta function of $F_{2,4}$ }}, in preparation.

\bibitem{Voll/03a}
C.~Voll, \emph{{Functional equations for local normal zeta functions of
  nilpotent groups}}, Geom. Funct. Anal., with an Appendix by A. Beauville, to
  appear (http://arxiv.org/abs/math.GR/0305362).

\end{thebibliography}

\end{document}